\magnification\magstep1
\baselineskip=18pt

\centerline{{\bf Remarks on complemented subspaces of von-Neumann
algebras}\footnote*{Supported in part by N.S.F. grant
DMS 9003550}} \vskip12pt \centerline {by  Gilles Pisier}
\vskip12pt \centerline {Texas A. and M. University and
Universit\'e Paris 6} \vskip12pt

{\bf Abstract}

In this note we include two remarks about
bounded ($\underline{not}$ necessarily contractive) linear projections on
a von Neumann-algebra. We show that if $M$ is a 
von Neumann-subalgebra of $B(H)$ which is complemented in
B(H) and isomorphic to $M \otimes M$ then $M$ is injective
(or equivalently $M$ is contractively complemented). We do
not know how to get rid of the second assumption on $M$.
In the second part,we show that any complemented 
reflexive subspace of a $C^*$- algebra is necessarily
linearly isomorphic to a Hilbert space.
\vfill\eject

1.
Let us start by defining the projection constant of a closed subspace  $M$ of a
Banach space  $B$  as 
$\lambda (M,B) =\inf(\|P\|,\quad P:B \rightarrow M)$
where the infimum runs over all bounded linear projections from $B$ into $M$.
The subspace $M$ is called complemented if there is a
bounded linear projection onto $M$, i.e. if $\lambda
(M,B)<\infty.$ ( We set $\lambda (M,B)=\infty$ when there
is no such projection.) Observe that if $M$ is weak*-closed
in $B$ then this infimum is attained, so that there exists
a bounded linear projection onto $M$ with norm 
$\underline{equal}$ to $\lambda (M,B)$. We might also
observe immediately that if $B=B(H)$ for some Hilbert
space $H$ and if  $M$ is a $C^*$ -subalgebra embedded
into  $B(H)$,then the projection constant $\lambda
(M,B(H)) $ does not depend on the embedding or on the
Hilbert space H, because all the $C^*$ morphisms  from M
into $B(H)$ extend with the
 same norm to operators defined on the whole of
$B(H)$,(cf.[A]). We can thus denote unambiguously
$\lambda(M)$ the projection constant in that particular
case.This is an invariant of the $C^*$-algebra M. \vskip4pt

Let us recall the definition of the projective tensor product $X\hat{ \otimes}Y$
of two Banach spaces $X$ and $Y$.For any $u$ in the linear
tensor product $X\otimes Y$ let

$$\|u\|_{\wedge}=\inf\Big\{\sum \|x_i\|\  \|y_i\| \Big\}$$
 
where the infimum runs over all the  representations of $u$ as a finite sum of
the form $u=\sum x_i \otimes y_i$.

We will use the following  well known lemma (which gives a dual
criterion for complementation):\vskip12pt

\noindent {\bf Lemma:}\quad Let $S$ be a weak*-closed subspace of a dual space
$X^*$ and let us denote by $S_{*}$ the predual $X/S^\perp$.Then the projection
constant  $\lambda (S)$ is equal to the smallest constant $\lambda$ such that 
for every $u$ in $S_* \otimes S$ we have
$$   |tr(u)| \leq \lambda \| u \|_{S_*  \hat{ \otimes} X^{*} }   . $$

We will make crucial use of the following proposition,  which could very well
have been observed by other researchers in Banach Space Theory.

\noindent {\bf Proposition1:}\quad Let $S^1$ and $S^2$ be
weak*-closed subspaces of two dual Banach spaces $X_1^*$
and $X_2^*$.We denote again by $S^1_*$ and $S^2_*$ the
respective preduals.Also, let $M$ be a weak*-closed
subspace of a dual Banach space $B$.We make the following
assumption (which means roughly that B is some kind of
tensor product of $X_1^*$ and $X_2^*$ and $M$
 some  tensor product of $S^1$ and $S^2$ ) :
There are a norm one operator $J:X_1^* \hat{\otimes}  X_2^* \rightarrow B$ and a
weak*-continuous operator $\phi:M \rightarrow (S^1_*\hat{\otimes}S^2_*)^* $
which is  also of norm one such that $J(S^1\otimes S^2) \subset M$
and the composition $\phi J$  coincides on  $S^1 \otimes
S^2$  with the canonical inclusion of  $S^1 \otimes S^2$
into $(S^1_*\hat{\otimes}S^2_*)^*$.\
Then $$\lambda(M,B) \geq \lambda(S^1,X_1^*) \lambda(S^2,X_2^*) .$$
Proof: \quad Fix $\epsilon >0$ and let $\lambda_1
=\lambda(S^1,X_1^*)-\epsilon$, and $\lambda_2 =\lambda(S^2,X_2^*)-\epsilon$.
By the above Lemma,  there are  $u_1$ and $u_2$  such that
$$u_1 \in S^1_*\otimes S^1,u_2 \in S^2_*\otimes S^2 $$
$$\| u_1\|_{S^1_*  \hat{ \otimes} X_1^{*} }=1 ,\quad
|tr(u_1)|> \lambda_1 ,\quad
\| u_2\|_{S^2_*  \hat{ \otimes} X_2^{*} }=1,\quad
|tr(u_2)|> \lambda_2 . $$
Then we use our assumption, we first denote by $\phi_*  :S^1_*  \hat{ \otimes}
S^2_*  \rightarrow M_*$ the preadjoint of $\phi$  and we note that its norm
is one . Let
 $$ u=(\phi_*\otimes J)(u_1\otimes u_2) \in M_*\otimes B .$$ 

Observe that, by our assumption,
$$u\in M_*\otimes M ,\quad \|u\|_{M_*  \hat{ \otimes}B } \leq 1 $$
and since $\phi J$ restricted to $S^1\otimes S^2$ is the canonical duality
mapping we also have $$tr(u)=tr(u_1) tr(u_2) >\lambda_1\lambda_2.$$Therefore by
the Lemma again we obtain $\lambda(M,B) > \lambda_1\lambda_2 $ which concludes
the proof.

The preceding proposition is tailor-made to apply to 
von Neumann algebras.Indeed, we can take for instance
$$X_1^*=B(H_1), \quad X_2^*=B(H_2) $$
where $H_1$ and $H_2$ are two Hilbert spaces, and we can
take for $S^1$ and $S^2$ two von Neumann subalgebras
respectively of $B(H_1)$ and $B(H_2)$.Then, if $H=H_1
\otimes H_2$ and if $M$ is the von Neumann tensor product
of $S^1$ and $S^2$ (i.e. the weak*-closure of $S^1 \otimes
S^2$ in $B(H)$ ) the above assumption is well known to be
verified (cf.e.g.[T] chapter 4.5, p.220).
Thus we have proved
 
\noindent {\bf Proposition 2:} \quad Let $M_1\subset
B(H_1),\quad M_2\subset B(H_2)$ be two von Neumann
subalgebras, and let $M=M_1\otimes M_2$ be their
von Neumann tensor product included in $B(H_1\otimes
H_2)$. Then
$$\lambda (M,B(H_1\otimes H_2)) \geq \lambda (M_1,B(H_1))
\lambda (M_2,B(H_2)).$$

By the remarks at the beginning of this section we
see that if $M$ is isomorphic (as a von Neumann algebra)
with $M\otimes M$ then $\lambda(M,B) $ must be equal
either to 1 or to infinity. Hence we obtain

\noindent{\bf Theorem1:}\quad Let $M$ be a von Neumann
subalgebra of $B(H)$. Assume $M$ isomorphic to $M\otimes
M$, (actually, it suffices to assume that $M\otimes
M$ embeds into
$M$ as a contractively complemented subalgebra). Then if
$M$ is complemented in $B(H)$, it must be contractively
complemented in $B(H)$.

The preceding theorem complements a well known result due
to Tomiyama [To] which says that when $M$ is contractively
complemented in $B(H)$, then $M$ is injective and the norm
one projection  is necessarily completely positive. It
would of course be nice to remove the extra condition, $M$
isomorphic to $M\otimes M$, but we do not see how to do
this.

As a corollary, we can deduce that there are
uncomplemented von Neumann subalgebras of $B(H)$. Indeed,
 let $N$ be any non injective von Neumann 
algebra  and let $M$ be  an infinite tensor product of
copies of $N$, then $M$ clearly is isomorphic to its
"tensor square" and is not injective since $N$ is
contractively complemented in $M$. Another approach
is to  define $N^n =N\otimes N\otimes ...\otimes N$ ($n$
times) and to define $M$ as the direct sum in the
$\ell_\infty$ sense of the sequence of spaces
$(N^n)_{n\geq 1}$. Then by Prop. 2, $\lambda (N^n) \geq
(\lambda (N) )^n$, hence if $\lambda (N) > 1$ (i.e. if
$N$ is not injective), we must have $\lambda (M) =\infty$.

Finally, we would like to mention that (although
everybody seems convinced that the answer is negative) we
could not find a counterexample to the following

\noindent {\bf Problem:}Let $M$ be a von Neumann
subalgebra of $B(H)$. Does every bounded operator from
$M$ into $B(H)$ extend to a bounded operator on the whole
of $B(H)$?

\vfill\eject
2.
 In the second part, we would like to record here an
 application of the main result of [P ],
namely the following (which answers a question raised by
G.Robertson)

\noindent {\bf Theorem 2:}
Every reflexive complemented subspace $S$ of a C*-algebra
$A$ is necessarily isomorphic to a Hilbert space.

Note that Hilbert spaces do appear as complemented
subspaces of $B(H)$. For instance the subspace of all
operators with values in a fixed one dimensional subspace
is obviously 1-complemented and isometric to $H$.

Proof:
Let $P:A\rightarrow S$ be a bounded linear projection
onto $S$. By a result of Jarchow [J] and by Theorem 3.2
(and the remark following it) in [P],
 we know that there are, 
a number $0<\theta <1$, a constant $C$, and a state $f$ on
$A$ such that  $$\forall x\in A,\quad \|P(x)\|\leq C\|P\|
(f(x^*x+xx^*))^{\theta/2}.\|x\|^{1-\theta}.$$

Since $\forall x \in S,\quad P(x)=x$, we deduce
$$\forall x \in S, \quad \|x\|\leq
{(C\|P\|)}^{1/\theta}.(f(x^*x+xx^*))^{1/2},$$
hence, finally
$$\forall x \in A, \quad 
(f({x^*x+xx^*}/2))^{1/2} \leq \|x\| \leq 2^{1/2}
C^{1/\theta}(f(x^*x+xx^*))^{1/2},$$
which shows that $S$ is indeed isomorphic to a Hilbert
space.

Note that, for the conclusion to hold, it suffices that
the subspace $S$ be of cotype $q$  where
$1/q=\theta/2$. If we denote by $C_q$ its cotype $q$
constant, we find the estimate
$$\|P\|\geq K(C_q)^{-1} d(S,H)^{2/q} ,$$
where $d(S,H)$ denotes the Banach-Mazur distance of $S$ to
a Hilbert space $H$ of the same dimension. It would be
interesting to know, for instance in the case when $S$ is
the $n$-dimensional space $\ell_q^n$ , what is the best
possible lower bound for the projection constant of $S$
inside a C*-algebra. The same question arises
naturally for the best constant of factorisation of the
identity of $\ell_q^n$ through a C*-algebra (or through
$B(H)$). It is not even clear that we can reduce these
questions to the single case of $B(H)$ .

\noindent {\bf Remark:}We would like to take this
opportunity to correct an error in [P]. In the middle
of page 123 of [P], it is erroneously stated that a
certain subset $N$ of a C*-algebra $A$ is a two sided ideal
of $A$, and that the quotient space $A/N$ is a
C*-algebra. This error does not affect the
rest of the paper [P]. Indeed, the subsequent discussion
makes perfectly good sense even if $A/N$ is merely a Banach
space. One more time, on page 124 of [P], a similar error
appears, when we state that $xf$ and $fx$ depend only on
the equivalence class of $x$ in $A/N$ , but here again the
error does not affect the argument because we only use the
fact that $2^{-1}(xf+fx)$ depends only on the equivalence
class of $x$ in $A/N$, and this is clearly correct.

\noindent {\bf Acknowledgement:} I am grateful to
G.Robertson for sending me his papers on his work related
to Theorem 2 above, and to G.Blower for stimulating
conversations on the same subject.

\vfill\eject

\centerline {\bf References}\vskip6pt
\item {[A]}  Arveson, W.B.:Subalgebras of
C*-algebras.Acta Math. 123 (1969) 141-224.
\vskip6pt
\item {[J]}  Jarchow, H.:Weakly compact operators on
C*-algebras.Math.Ann.273 (1986) 341-343.
\vskip6pt
\item {[P]} Pisier, G.:Factorization of operators through
$L_{p\infty}$ or $L_{p1}$ and non-commutative
ge-neralizations.Math.Ann.276 (1986) 105-136.
\vskip6pt
\item {[T]} Takesaki, M.:Theory of 0perator Algebras,
vol.1.Springer Verlag.(1979)

 \item {[To]} Tomiyama, J.:On the projection of norm
one in W*-algebras.Proc.Japan Acad. 33 (1957) 608-612.

\vskip24pt
Math. Dept.

Texas A. and M. University

College Station, TX 77843, USA

and

Equipe d'Analyse

Universit\'e Paris 6

Tour 46, 4\`eme \'etage

75230 Paris Cedex 05, FRANCE

 \bye